**Do Machine Learning–Guided Proactive Full Capacity Protocols Mitigate Emergency Department Overcrowding? A Simulation-Based Evaluation**

Khalid Aram[1], Orhun Vural[2], Bunyamin Ozaydin[3,4], James Booth[5], Abdulaziz Ahmed[3,4*]

[1]School of Business & Technology, Emporia State University, Emporia, KS, United States
[2]Department of Electrical and Computer Engineering, University of Alabama at Birmingham, Birmingham, AL, United States
[3]Department of Health Services Administration, School of Health Professions, University of Alabama at Birmingham, Birmingham, AL, United States
[4]Department of Biomedical Informatics and Data Science, Heersink School of Medicine, University of Alabama at Birmingham, Birmingham, AL, United States
[5]Department of Emergency Medicine, University of Alabama at Birmingham, Birmingham, AL, United States

## Abstract

*Objective:* Emergency department (ED) overcrowding threatens care quality, yet Full Capacity Protocols (FCPs) are typically activated only after crowding thresholds have been exceeded. This study evaluated whether machine learning–guided proactive FCP activation improves ED crowding control compared with conventional reactive activation.

*Materials and Methods:* We developed and validated an integrated machine learning–discrete-event simulation framework representing ED patient flow and capacity management. An embedded machine learning model predicted hourly waiting-room crowding six hours ahead and triggered proactive FCP activation. Using empirically calibrated operational parameters from a large U.S. hospital, we compared proactive and reactive deployment strategies over a 365-day simulation horizon with 10 independent replications. Outcomes included ED length of stay, wait time, percentage of time above the crowding threshold, number of FCP activations, and protocol duration.

*Results:* Proactive deployment produced statistically significant reductions in the percentage of time above the crowding threshold and average protocol duration across all threshold settings. It also reduced crowding variability and extreme congestion events, resulting in smoother ED operations. Although proactive deployment generally improved ED length of stay and wait time, these differences were not statistically significant.

*Conclusion:* Embedding machine learning predictions within an operational simulation enables proactive, prediction-guided capacity management before severe congestion develops. These findings support the feasibility and potential value of proactive FCP deployment as a decision-support strategy for improving ED crowding control and provide evidence to support prospective implementation in clinical practice.

[*] **Corresponding Author,** Abdulaziz Ahmed, PhD, Department of Health Services Administration, School of Health Professions, University of Alabama at Birmingham, Birmingham, AL 35233, United States, Email: aahmed2@uab.edu

**Key words:** Emergency Overcrowding, full capacity protocol, Simulation, Machine Learning.

## Lay Summary

Emergency departments often become overcrowded, leading to longer waits and delays in patient care. Hospitals typically respond only after overcrowding has already occurred. In this study, we developed an artificial intelligence (AI) system that predicts when overcrowding is likely to happen several hours in advance. We combined this prediction system with a computer simulation of a hospital emergency department to test whether acting earlier would improve operations. Our results showed that using AI to activate full capacity protocols before severe crowding develops reduced the amount of time the emergency department remained overcrowded and shortened the duration of emergency capacity measures. These findings suggest that AI can help hospitals respond more proactively to rising demand and may improve emergency department operations without disrupting patient care.

## Introduction

Emergency Department (ED) overcrowding is a critical issue in healthcare systems worldwide, leading to significant negative impacts on patient care and hospital operations. According to the most recent data from the Centers for Disease Control and Prevention, there were approximately 140 million ED visits, equating to about 43 visits per 100 persons.[1] Of these visits, 18.3 million (13.1%) led to hospital admission, and 2.8 million (2%) ED visits led to admission to a critical care unit.[2] EDs are overcrowded due to many factors, including the complexity of complaints and injuries,[3,4] resource limitations,[5] and large patient volumes.[6] ED overcrowding is associated with poor healthcare outcomes. It leads to delays in diagnosis and treatment, which can result in worse clinical outcomes, increased morbidity, and greater illness severity.[7,8] For example, patients with acute coronary syndrome who presented to overcrowded EDs had a significantly higher rate of serious complications than those who presented to non-crowded ED (6% vs 3%), including death, late myocardial infarction, and hypotension.[9] Sprivulis et al.[10] observed a significant linear relationship between ED overcrowding and patient mortality based on three years of data from three large hospital systems. Additionally, ED overcrowding can also negatively impact mental health, as chaotic environments can worsen symptoms for patients with mental health conditions. These patients may also experience longer wait times and boarding due to limited psychiatric resources, further affecting their well-being.[11,12] Beyond mortality, ED overcrowding is a major obstacle to achieving national quality measure goals and affects the national capacity for disaster response.[13] It causes stress for medical staff,[14,15]increasing the likelihood of medical errors.[16] Delays in treatment and long ED wait times negatively affect patient satisfaction and quality of care.[17,18]

ED overcrowding is not only multifactorial, but also a complex, multi-level patient-flow problem. At the hospital level, overcrowding is driven by inefficient patient flow, limited resources such as inpatient beds, and poor coordination between the ED and inpatient units. The Emergency Medicine Practice Committee of the American College of Emergency Physicians (ACEP) reported that ED overcrowding is a hospital-wide patient-flow problem rather than an isolated ED problem.[19,20] Many studies also suggest that ED crowding stems from bottlenecks caused by inefficient patient flow, particularly boarding of admitted patients in the ED.[21-24] In addition, the Agency for Healthcare Research and Quality (AHRQ)'s guide for hospitals on improving patient flow and reducing ED crowding,[25] along with reports from the Institute of Medicine, the Institute for Healthcare Improvement, and the Joint Commission, emphasizes patient-flow improvement as a key strategy for mitigating ED overcrowding. [26-28,20,29] A key approach suggested by ACEP is the Full Capacity Protocol (FCP), an internationally recognized hospital-wide escalation process used to coordinate actions between the ED and inpatient units.[20,29] FCPs typically include multiple levels of intervention triggered by patient flow measures (PFMs), such as the number of boarded patients, waiting-room counts, or other indicators of ED crowding. However, current FCP practice is largely reactive because activation decisions rely on near real-time PFM values. As a result, FCP interventions are often initiated only after crowding has already developed. This timing is a major limitation because many FCP actions, such as coordinating staffing, creating hallway treatment spaces, expediting discharges, or restricting non-urgent transfers, require preparation time. Therefore, patient-flow teams often must prepare and respond while the ED is already overcrowded.

Predictive analytics offers a way to address this timing limitation by forecasting the PFMs that guide FCP activation. These forecasts can add anticipatory capability to FCP, transforming it from a reactive tool into a proactive strategy that supports earlier, better-coordinated responses before conditions deteriorate. By accessing predicted PFM values in advance, the patient flow coordination team can prepare and implement interventions before overcrowding has fully developed. Recent artificial intelligence (AI) and machine learning (ML) studies have demonstrated the ability to forecast ED waiting room count, length of stay, bed occupancy, and operational bottlenecks using a wide range of time-series and deep learning models.[30-33] Short-horizon forecasting is especially relevant for operational decision-making, as many crowding mitigation strategies, including hallway-bed deployment, staffing adjustments, and transfer restrictions, have activation lead times measured in hours, not minutes. However, most ML studies evaluate only predictive performance and do not test whether acting on these predictions improves downstream system behavior, leaving an important gap between prediction and measurable operational impact.

Simulation models can be used for rigorous evaluation of how predictive models influence ED operations once embedded in the ED workflow. Discrete-event simulation (DES) is widely used in ED operations studies because it captures patient flow dynamics, resource constraints, queueing behavior, and the interaction between ED and inpatient units.[34,35] Prior DES studies have tested interventions such as physician-in-triage, fast-track units, bedside registration, early discharge strategies, and dynamic bed allocation.[36] These studies emphasize the importance of integrated hospital and ED modeling as they consistently show that ED performance depends on both micro-level processes (triage efficiency, treatment times) and macro-level constraints (inpatient occupancy, discharge timing). Yet attempts at embedding a real-time, ML-based prediction model directly inside the simulation loop remain limited. Our study advances the growing class of hybrid decision-support frameworks that integrate ML and simulation to evaluate operational policies. The literature includes limited simulation-based studies of prediction-driven interventions. For instance, Mišić et al.[16] embedded a readmission-risk model into a hospital workflow simulator to assess its real-world clinical and operational impact. Other methodological precedents include the integration of admission-prediction classifiers directly inside patient-flow DES models to dynamically influence routing and system behavior during simulation.[17] Building on this paradigm, our study advances ED crowding research by integrating a short-horizon ED crowding prediction model into a validated ED simulation model. The forecasts are used to trigger proactive FCP deployment before crowding thresholds are exceeded. This predictive-control design allows us to evaluate how earlier action affects patient-flow dynamics, crowding duration, and overall ED stability.

The primary objective of this study is to evaluate whether earlier, prediction-informed activation of the FCP can improve ED operations. In other words, our question is: *does the timing of FCP deployment matter?* To address this question, we developed a DES model of a large ED at an academic medical center in the southeastern United States. The model integrates an ML-based hourly crowding prediction model directly into the simulation cycle, enabling real-time, proactive activation of FCP. Our simulation model incorporates patient arrival distributions, ESI acuity distributions, treatment and triage times, real staffing schedules, inpatient bed capacities, and LOS distributions across multiple units. Unlike traditional reactive activation, which occurs only once thresholds are exceeded, the proactive strategy triggers FCP deployment based on predicted overcrowding six hours ahead. We compare both strategies across multiple replications over a full simulated year and evaluate differences in waiting-room dynamics, volatility, peak crowding events, boarding burden, holds, wait times, and the overall ED length of stay.

Our results show that prediction-guided proactive FCP activation reduces extreme crowding peaks, increases system stability, and shortens crowding duration compared with conventional reactive deployment. These findings demonstrate that ML-guided capacity management, when operationalized through simulation and embedded decision logic, can strengthen ED resilience and reduce instability in tightly coupled ED–inpatient systems. More broadly, our study provides a methodological framework for translating predictive analytics into operational improvements using simulation-based evaluation.

## Methods

We designed a structured framework, illustrated in Figure 1, for our methodology that integrates three main parts: data processing, simulation-based pilot testing, and predictive model development. The framework was designed to evaluate the effectiveness of our proposed proactive FCP decision-support system. Operational data from the ED Encounters and ED Flow information were first acquired and processed to derive probability distributions for arrivals, care-process times, acuity levels, dispositions, inpatient length of stay (LOS), and inpatient unit assignments (Appendix Tables A2–A7). These distributions formed the input parameters for the discrete-event simulation model, which replicates patient flow and resource constraints within the ED and inpatient processes. The simulation environment then served as a pilot-testing platform for comparing reactive and prediction-guided proactive FCP deployments. We ran the simulation model without FCP first to generate hourly wait count data. We then used the simulated data to train and test a machine learning model that was subsequently used as a part of the FCP control module within the simulation model to predict near-term ED waiting counts. The ML model and simulation were integrated so that the predictive engine could drive proactive activation decisions within the simulated ED. Our framework provides a safe, data-driven means to evaluate the proposed decision-support system under realistic operational conditions prior to real-world implementation of proactive FCP. The next part of the methods section focuses on an in-depth description of the different parts of the framework.

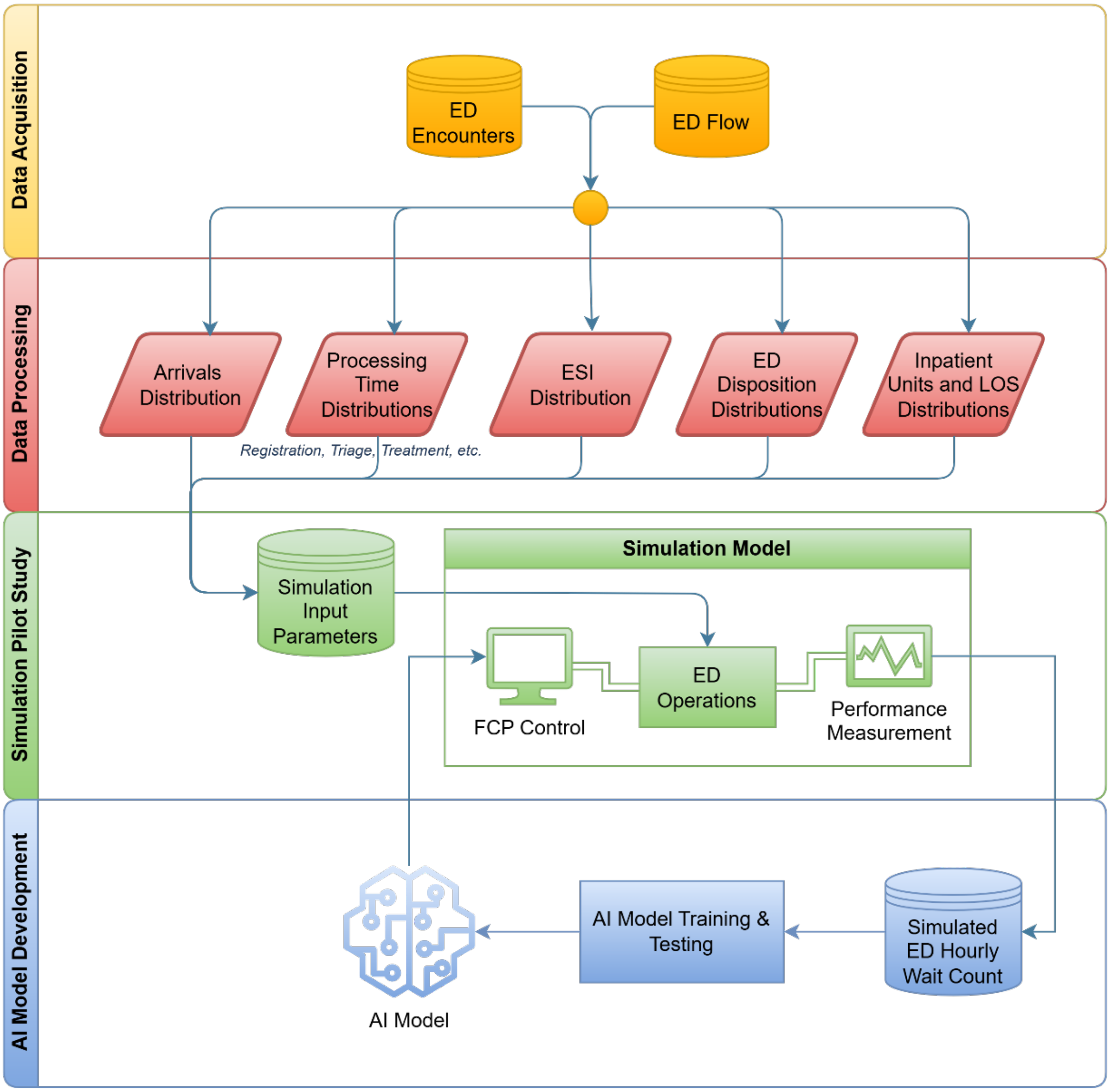

Figure 1: Research framework.

## Data Acquisition and Processing

We generated parameters for our simulation model using ED encounters data from our partner hospital between January 2019 and July 2023, totaling 308,196 patient visits. We excluded data from the COVID-19 period. We also collected data on hospital inpatient units and their capacity, inpatient LOS, and admission and discharge records. We then analyzed the data and obtained parameters for the probability distributions that best fit the different data elements such as arrivals, registration time, triage time, treatment time, ESI distribution, ED nurse staffing schedule, admission distribution, inpatient unit distribution by ESI, and inpatient LOS. These distributions are included in the appendix and in the published project code of the simulation model on GitHub.[39]

## Crowding Predictive Model

We developed a machine learning model to predict hourly ED wait counts six hours ahead based on the most recent 48 hours of hourly waiting count data. Model inputs, prediction frequency, and FCP trigger rules are summarized in Appendix Table A8. We trained the model using simulated hourly wait counts after excluding output of the simulation 72-hour warmup period. We explored multiple time-series forecasting models using five-fold cross validation and reported average and standard deviation of predictive performance. The model comparison shown in Table 1 shows that Support Vector Regression (SVR) with a Radial Basis Function (RBF) kernel provides the lowest average Mean Absolute Error (MAE) of 2.46 patients and the highest average $R^2$ of 83.3%.

SVR offers relatively lower computational cost at inference compared with more complex models such as neural networks. This contributes to substantial reduction of the computational burden as the predictive model is embedded within the simulation and executed repeatedly during each run.

Table 1: Performance comparison of wait count predictive models.

| Model | AVG MAE | MAE STDEV | $R^2$ AVG | $R^2$ STEV |
|---|---|---|---|---|
| SVR RBF | 2.46 | 0.266 | 0.833 | 0.0264 |
| Linear Regression | 2.62 | 0.218 | 0.815 | 0.0164 |
| Ridge Regression | 2.62 | 0.219 | 0.815 | 0.0165 |
| Gradient Boosting | 2.62 | 0.252 | 0.811 | 0.0246 |
| Random Forest | 2.64 | 0.278 | 0.806 | 0.0264 |

## Simulation Model

This section provides a description of the conceptual logic of our simulation model in addition to details regarding the modeling approach. We applied Discrete-event Simulation (DES), an approach well-suited for process-based simulation applications.[36] The goal of the model is to evaluate the effectiveness of our proposed proactive FCP in comparison to the conventional reactive approach. The flowchart in Figure 2 illustrates the conceptual model. The model includes three main portions: ED front-end operations, ED-back-end operations, and FCP control logic. Core simulation settings, ED capacity assumptions, staffing schedules, and FCP operating assumptions are summarized in Appendix Table A1. We consider the level of detail in our model to be moderate, which is suitable for the scenario analysis task at hand. We implemented our simulation model in Python, following an event-scheduling approach, where the model maintains a calendar of events and advances the simulation time to the time of the nearest future event on the calendar. The simulation code can be found on GitHub.[39]

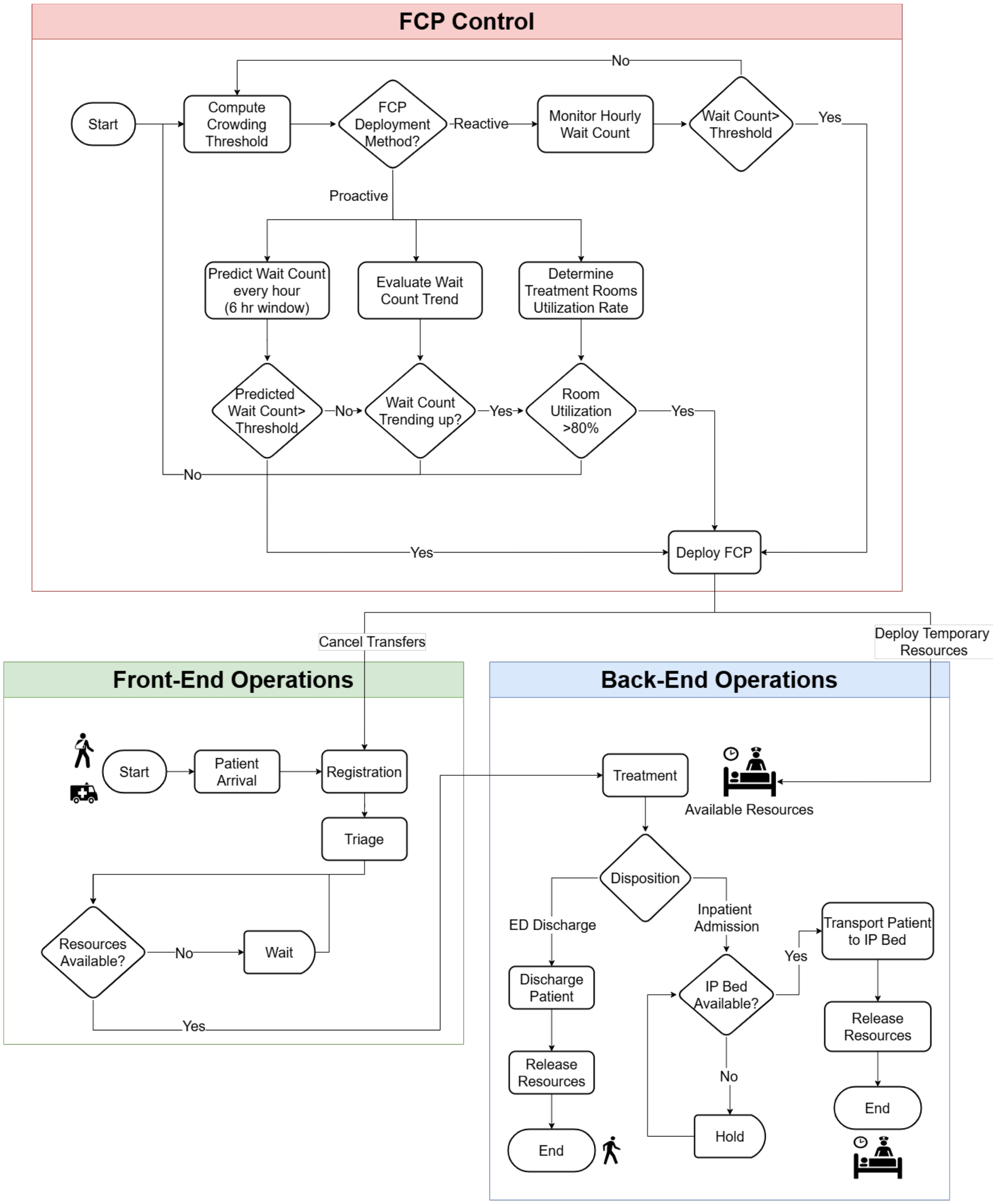

Figure 2: Conceptual ED simulation model.

The FCP control logic in Figure 2 illustrates how the FCP is activated. The prediction inputs and reactive and proactive FCP trigger rules are summarized in Appendix Table A8. The ED crowding threshold is computed dynamically throughout the simulation based on the hourly wait counts from the most recent 48 simulation hours. Specifically, the crowding threshold is calculated as $\mu + z\sigma$, where $\mu$ is the average hourly wait count, $\sigma$ is the standard deviation of hourly wait count, and z is a parameter that sets an upper limit (worst-case scenario) on the threshold. The FCP control module applies a different FCP activation logic depending on the preset firing mode: reactive or proactive. The reactive mode activates the FCP when the wait count reaches or exceeds the current crowding threshold. In proactive mode, FCP is activated if one of two conditions is met. The first condition is if the 6-hour predicted wait count is larger than the threshold for two consecutive hours. The second condition is if the real-time hourly wait count exhibits a positive trend for four hours AND the current treatment room utilization rate exceeds 80%.

**FCP Design**

The primary objective of this study was to evaluate whether earlier, prediction-informed activation of the Full Capacity Protocol (FCP) can improve ED operations, rather than to redesign the protocol itself. In essence, our question was: *does the timing of FCP deployment matter?* To isolate the effect of timing, we held the FCP intervention bundle constant across all scenarios and modified only the activation logic, comparing the traditional reactive approach with our proposed proactive strategy triggered by anticipated crowding. The FCP design in our simulation includes three interventions: deploying ten hallway beds, deploying five temporary nurses, and cancelling inbound transfers (3% of total ED visits) (Appendix Table A1). The following is a list of FCP timing parameters:

- FCP activation delay = 3 hours (time needed to activate FCP)
- FCP minimum active duration = 6 hours
- FCP cooldown duration = 6 hours (minimum time between FCP activations)
- FCP deactivation condition: 60% of activation threshold achieved for 4 consecutive hours.

**Experimental Setup**

We evaluated the FCP deployment scenarios using our DES model of ED operations run over a full 365-day horizon, allowing the model to capture long-range crowding patterns and high-intensity crowding events. We replicated each scenario 10 times using different random seeds and calculated average and standard deviation of all performance measures. The simulation included a warm-up period of 72 hours, during which system conditions were allowed to stabilize. All key performance measures (e.g., wait counts, holds, ED length of stay) were collected only after the warm-up window.

To evaluate the robustness of the proposed proactive FCP strategy, we conducted a targeted sensitivity analysis on key protocol parameters governing the timing, frequency, duration, and intensity of capacity interventions. This analysis assesses whether the benefits of prediction-guided deployment persist across alternative configurations and identifies conditions under which these effects strengthen or weaken.

**Statistical Analysis**

Each simulation scenario was run for 10 independent replications using different random seeds. Outcomes were summarized as the mean and standard deviation across replications. Figures report mean values with 95% confidence intervals calculated across the 10 replications. For model validation, simulated and observed ED LOS were compared using a two-sided two-sample t-test, with statistical significance assessed at $\alpha = 0.05$. For comparisons between reactive and proactive

FCP strategies, differences in simulation outcomes were assessed using two-sided t-tests across replications, with $\alpha = 0.05$. Actual P values were reported where statistical significance was evaluated.

# Results

## Simulation Model Validation

Model validation is a critical step in the development of simulation models. Validation helps in assessing to what degree the simulation is close to the real system being studied. To validate our simulation model, we compared the ED length of stay (LOS) times, i.e., ED length of stay, produced by the simulation model to the actual ED LOS data. We compared the ED LOS data from one year of actual data to one year of simulated data. Visual inspection indicates that the simulated LOS distribution closely replicated the shape and dispersion of the observed LOS distribution as shown in figure 3. We also conducted a t-test of the difference in mean ED LOS between simulated and observed data. We used a two-sample t-test assuming equal variances. The p-value from the test is 0.4848, indicating no statistically significant difference between actual and simulated ED LOS means. The mean LOS error between simulated and observed data is 0.21% (330.0 vs. 329.3 minutes). Figure 3 further shows the similarity of actual and simulated ED LOS values through the shapes of their histograms.

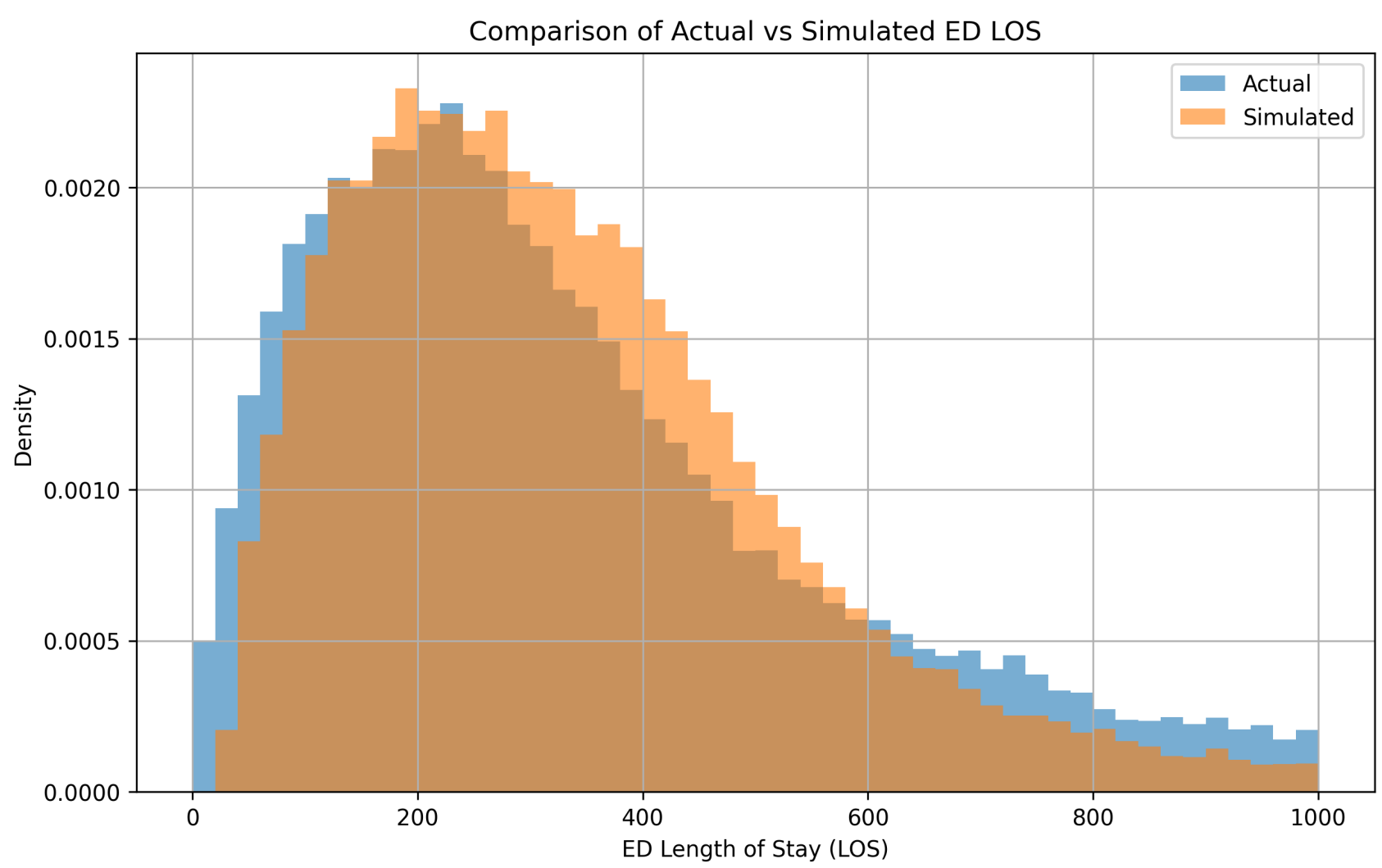

Figure 3: Histograms of actual and simulated ED LOS values.

## Simulation Results

In this section, we present simulation results in terms of a comparison between *reactive* (base-case) and *proactive* FCPs. We ran our model at three different levels of the ED crowding threshold. Crowding threshold is one of the critical parameters used to control FCP deployment. It is calculated dynamically throughout the simulation based on the observed recent waiting counts. These levels are $\mu + 1\sigma$, $\mu + 2\sigma$, $\mu + 3\sigma$, where $\mu$ and $\sigma$ are the average and standard deviation of the ED hourly wait count for the most recent 48 hours of simulation time, respectively. The $\sigma$ level controls the sensitivity of FCP activation. A lower standard deviation $\sigma$ produces lower activation threshold, which contribute to more frequent FCP activations. Table 2 includes definitions of the performance measures calculated by the simulation model which include ED operational measures and FCP frequency and duration measures. Figures 4–6 show the simulation results comparing reactive and proactive FCP deployment across the three crowding-threshold settings, using confidence interval plots of key ED performance measures.

Table 2: Definitions of ED simulation performance measures.

| Measure | Definition |
|---|---|
| AVG ED LOS | *Average total time from patient registration to final disposition (discharge or inpatient admission).* |
| AVG Wait Time | *Average delay between patient registration and placement into a treatment room.* |
| % Waits > Threshold | *Percentage of hourly wait counts exceeding the crowding threshold ($\mu + z\sigma$).* |
| Total Deploys | *Number of FCP activations during the simulation run.* |
| AVG Protocol Duration | *Average length of time FCP remains active once initiated.* |

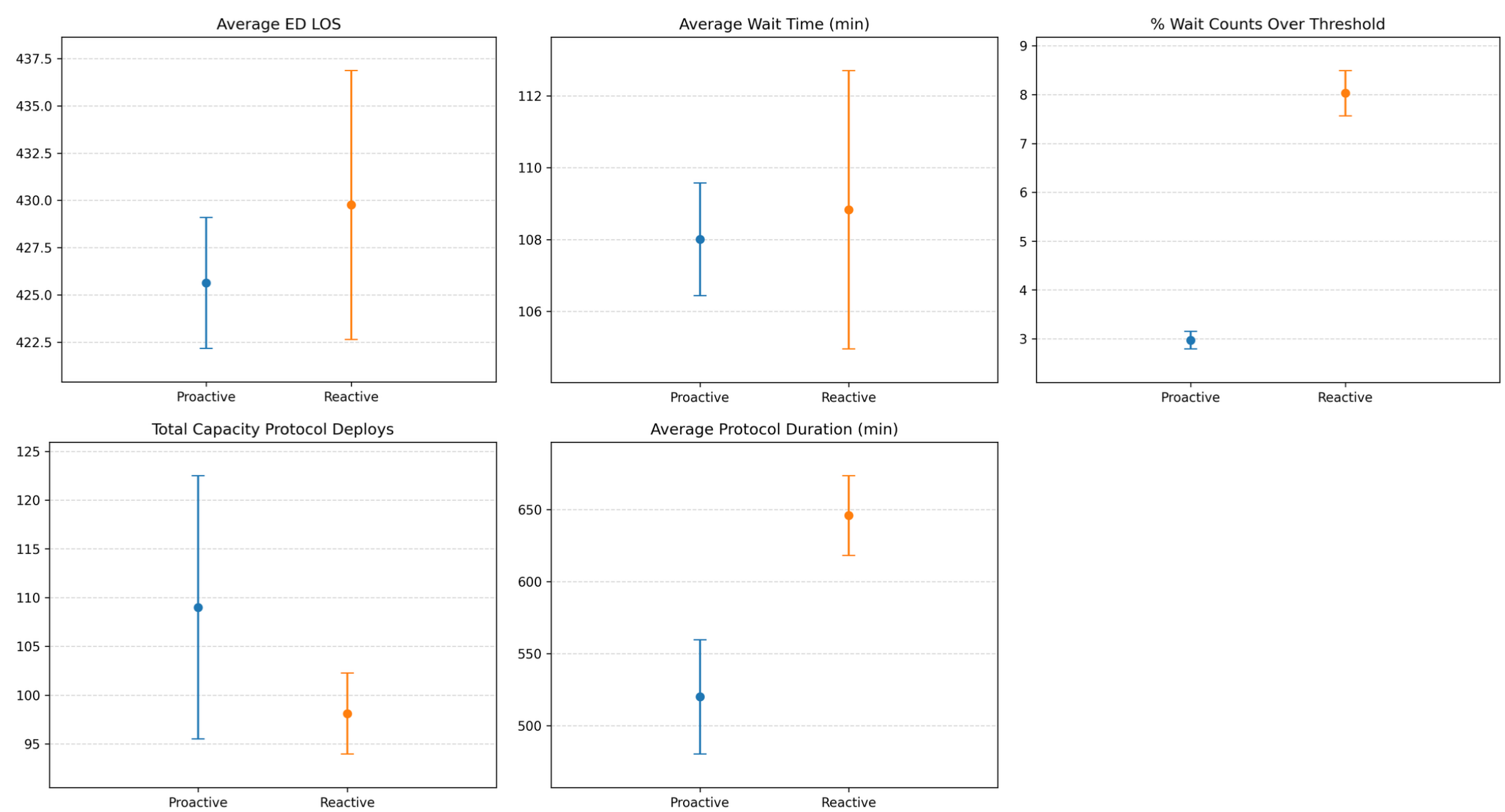

Figure 4: 95% confidence interval plots of simulation performance measures for reactive and proactive FCPs. (crowding threshold = $\mu + 2\sigma$)

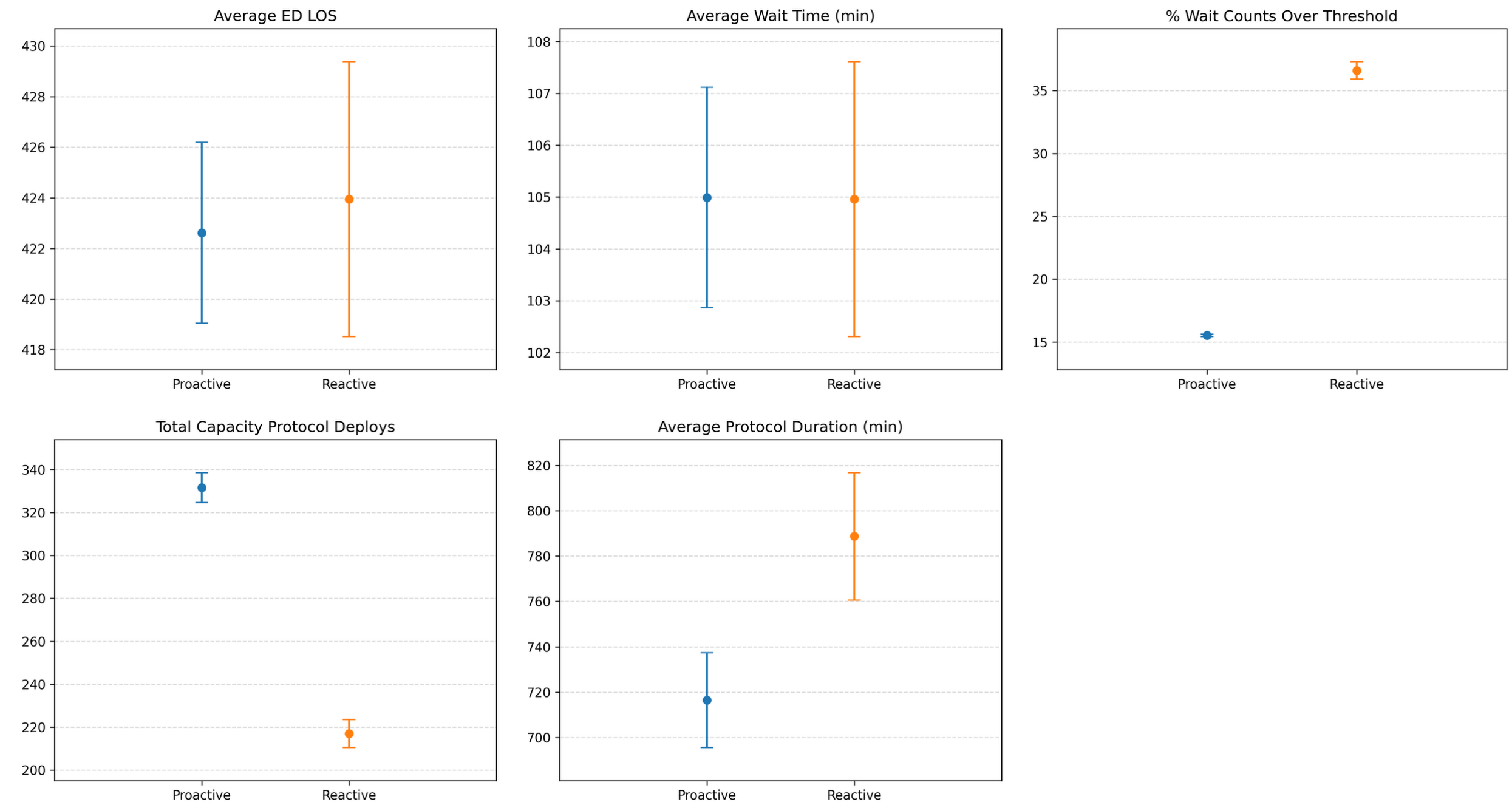

Figure 5: 95% confidence interval plots of simulation performance measures for reactive and proactive FCPs. (crowding threshold = $\mu + 1\sigma$)

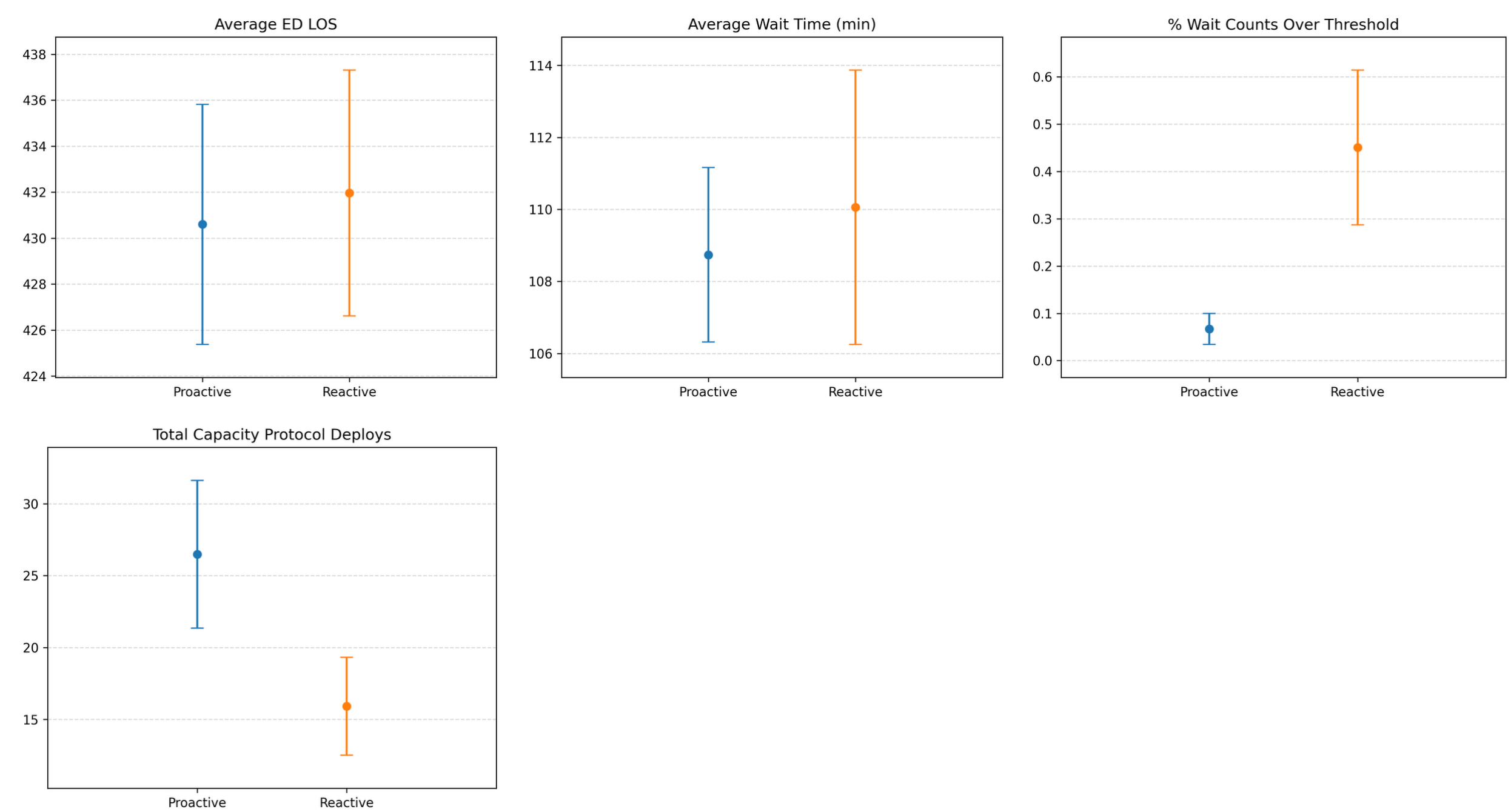

Figure 6: 95% confidence interval plots of simulation performance measures for reactive and proactive FCPs. (crowding threshold = $\mu + 3\sigma$)

At the moderate crowding threshold of $\mu + 2\sigma$, where protocol activation is triggered when the hourly waiting-room count exceeds the recent 48-hour average plus two standard deviations, proactive FCP deployment demonstrated measurable differences in operational outcomes over the reactive strategy (Figure 4). Proactive activation reduced average ED LOS by ~4 minutes and reduced average hourly wait time by close to one minute. The most substantial difference was observed in the percentage of time spent above the crowding threshold, which was markedly lower under proactive deployment (~5% reduction), suggesting better control of severe crowding situations. Although proactive deployment triggered the protocol around 10 times more than reactive deployment on average, activations were ~125 minutes shorter in duration. This suggests that earlier activation may reduce protocol duration by responding to crowding before it becomes prolonged. These findings suggest that prediction-guided proactive deployment yields smoother operations and more effective crowding control than the traditional reactive approach.

At the more sensitive crowding threshold of $\mu + 1\sigma$, where the FCP is activated when the hourly waiting-room count exceeds the recent 48-hour average by one standard deviation, proactive deployment produced operational effects that were stronger than those observed at the $\mu + 2\sigma$ setting (Figure 5). Because this lower threshold triggers capacity interventions earlier and more frequently, the proactive strategy reduced the percentage of time the ED operated above the crowding threshold to a much greater extent than in the $2\sigma$ scenario (~21% reduction). As expected, proactive deployment at this sensitive threshold resulted in a higher number of protocol activations than both the reactive strategy and the $2\sigma$ setting, yet the average protocol duration remained shorter.

At the least sensitive crowding threshold of $\mu + 3\sigma$, where the FCP is activated only when the hourly waiting-room count exceeds the recent 48-hour average by three standard deviations, the operational differences between proactive and reactive deployment were smaller than those observed at the $\mu + 2\sigma$ and $\mu + 1\sigma$ settings (Figure 6). Because this higher threshold allows the ED to tolerate substantially more crowding before any intervention is triggered, both strategies activated the protocol infrequently. However, even in this conservative setting, proactive deployment still resulted in fewer hours spent above the crowding threshold and slightly shorter protocol durations. Overall, the $\mu + 3\sigma$ results show that when the activation rule is highly restrictive, the operational impact of proactive decision support diminishes, reinforcing that the value of prediction-guided deployment is greatest at more responsive thresholds.

Across all three crowding-threshold settings, proactive FCP deployment consistently produced statistically significant reductions in the percentage of time spent above the crowding threshold and in average protocol duration. The average ED LOS and wait times were lower under proactive FCP, but these differences were not statistically significant. Hence, the main benefit of improving the timeliness of activation through proactive FCP lies in reducing severe crowding exposure and shortening protocol duration, rather than producing statistically significant reductions in average LOS or wait time.

**Sensitivity Analysis Results**

To assess the robustness of the proposed proactive FCP deployment strategy, we conducted a targeted sensitivity analysis on key protocol design parameters that control the timing, frequency, duration, and intensity of capacity interventions. This analysis also allows us to determine whether the observed benefits of proactive, prediction-guided FCP deployment remain consistent under alternative operational configurations and to identify conditions under which these benefits may strengthen or diminish. We evaluated the effects of variations in activation delay, cooldown duration, minimum active duration, resource levels (hallway beds and temporary nurses), and the strictness of FCP deactivation rules. Activation delay represents the time required to implement FCP after a trigger occurs, cooldown duration represents the minimum time required before the protocol can be reactivated after deactivation, minimum active duration represents the shortest time the protocol must remain active once initiated, and deactivation rule strictness defines the conditions required to end the protocol. We developed simulation scenarios of low and high levels of each parameter relative to the baseline configuration while maintaining all other settings constant. Table 3 shows a summary of the sensitivity analysis scenarios. All sensitivity scenarios were run under the $\mu + 2\sigma$ crowding threshold**.** In addition, selected interaction scenarios were tested to explore combined effects between parameters that may jointly influence the system response. Performance was evaluated using the same set of ED and FCP operational metrics introduced earlier, and comparisons between proactive and reactive deployment were preserved across all scenarios.

Table 3: Sensitivity analysis scenarios.

| **Parameter** | **Baseline** | **Low Level** | **High Level** |
|---|---|---|---|
| Activation Delay | 3 h | 1 h (Fast) | 5 h (Slow) |
| Cooldown Duration | 6 h | 3 h (Short) | 9 h (Long) |
| Minimum Active Duration | 6 h | 3 h (Short) | 9 h (Long) |
| Additional Resources (Bundle) | 10 beds + 5 nurses | 5 beds + 2 nurses | 15 beds + 8 nurses |
| Deactivation Rule Strictness | 60% of threshold factor for 4 hours | 70% of threshold for 3 hours (easier deactivation) | 50% of threshold for 5 hours (harder deactivation) |

Figure 7 shows the ED performance results under low, baseline, and high settings for key FCP parameters. The plots show that the effectiveness of proactive FCP deployment is sensitive to key protocol design parameters, particularly those controlling timing, resource availability, and deactivation rules. Across the sensitivity scenarios, proactive deployment consistently maintains a lower percentage of time above the crowding threshold compared with reactive deployment. This was the clearest and most consistent benefit of the proactive strategy. In contrast, differences in average ED LOS and average wait time were smaller and less consistent across parameter settings. The results suggest that higher resource levels tend to strengthen the benefits of proactive deployment in terms of ED LOS and average wait times, while reactive deployment did not show consistent improvement in these measures as resource levels increased. This suggests that when the ED is already overcrowded, additional resources alone may not be sufficient to improve outcomes once crowding has developed. On the other hand, proactive intervention may reduce severe crowding and reduce the resources needed to mitigate anticipated crowding. These findings suggest that proactive FCP primarily improves control of severe crowding events rather than uniformly reducing average throughput measures. Overall, these results highlight the importance of aligning prediction-driven decisions with operational feasibility, including activation delay, cooldown duration, resource levels, and deactivation criteria.

ED LOS | Average Wait Time | % Wait over Threshold

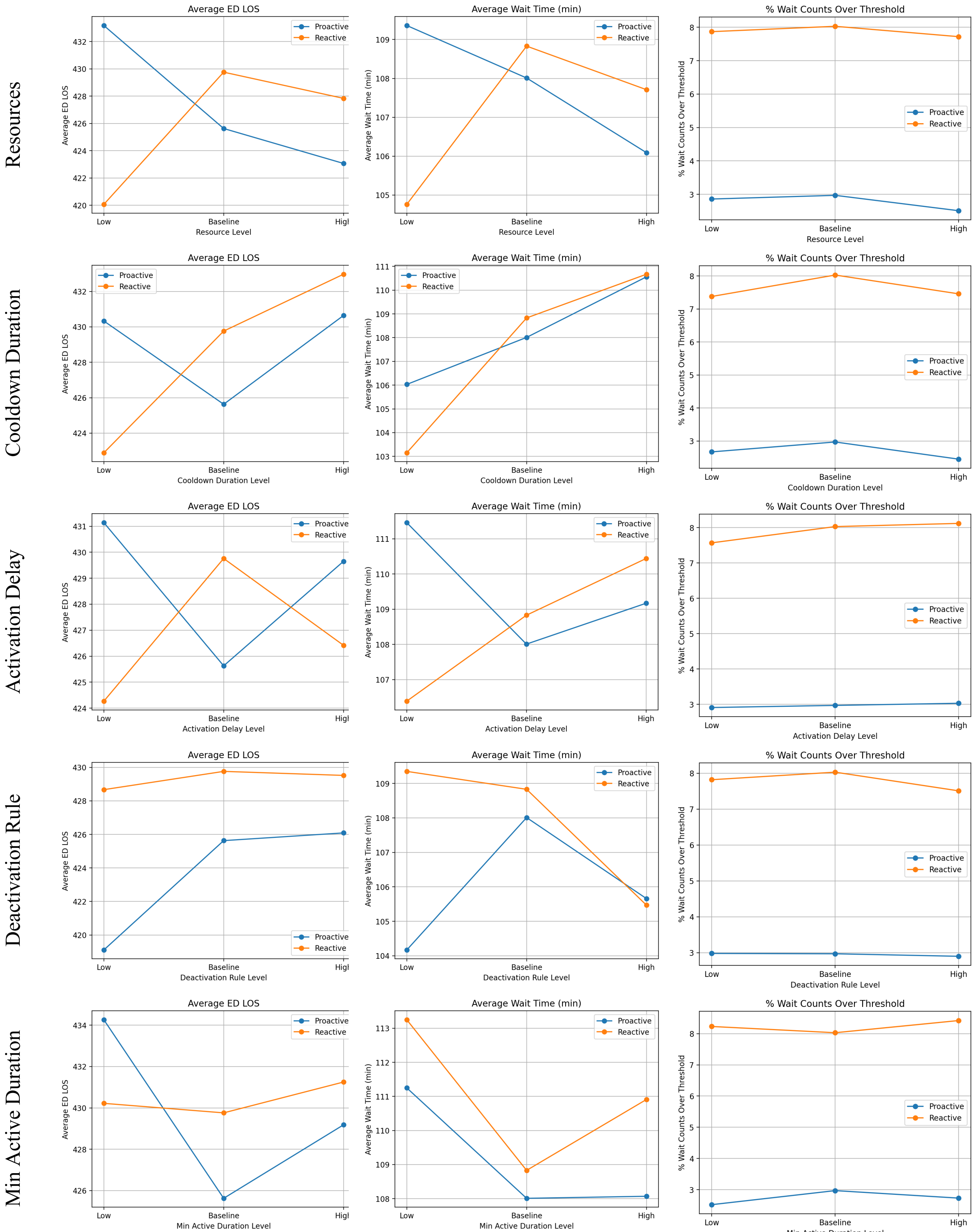

Figure 7: ED performance results from the sensitivity analysis of FCP parameters.

Figure 8 shows the results of FCP deployment metrics across the tested FCP parameter configurations. The plots highlight systematic differences in how proactive and reactive strategies utilize capacity interventions. Proactive deployment generally results in a higher number of protocol activations, particularly under more aggressive configurations such as shorter cooldown periods or shorter minimum durations. Despite this increased frequency, the average duration of each deployment is generally shorter than the reactive approach, which highlights the importance of targeted and timely interventions. In contrast, reactive deployment tends to rely on fewer but longer activations because FCP activation occurs only after crowding has already exceeded the threshold. These patterns are consistent across most parameter variations and suggest that proactive control shifts the system from infrequent, prolonged interventions toward more frequent, shorter adjustments, which may help limit prolonged crowding exposure. However, the duration advantage of proactive deployment was not uniform across all settings, particularly when activation delay increased, indicating that deployment timing remains an important determinant of protocol efficiency.

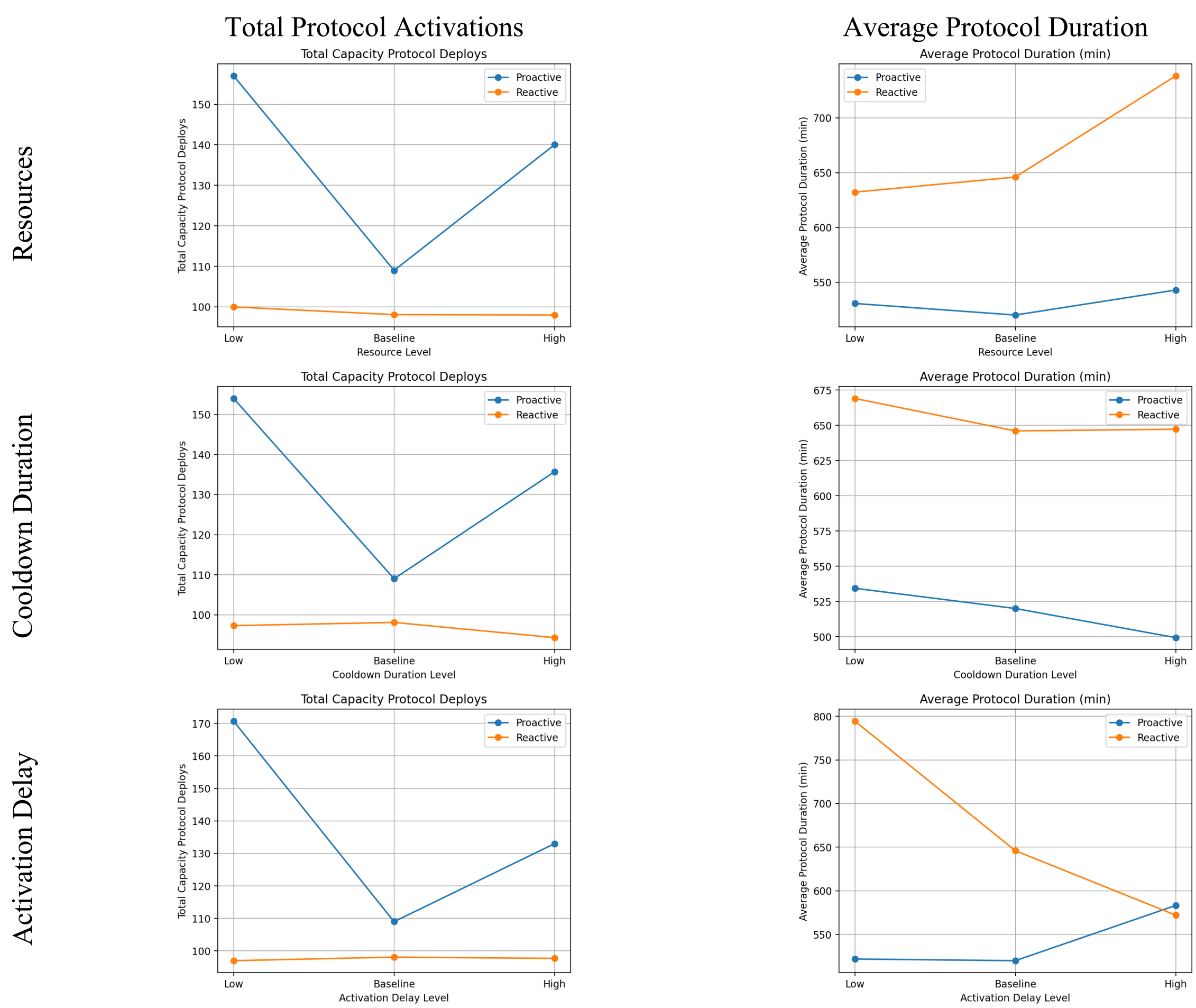

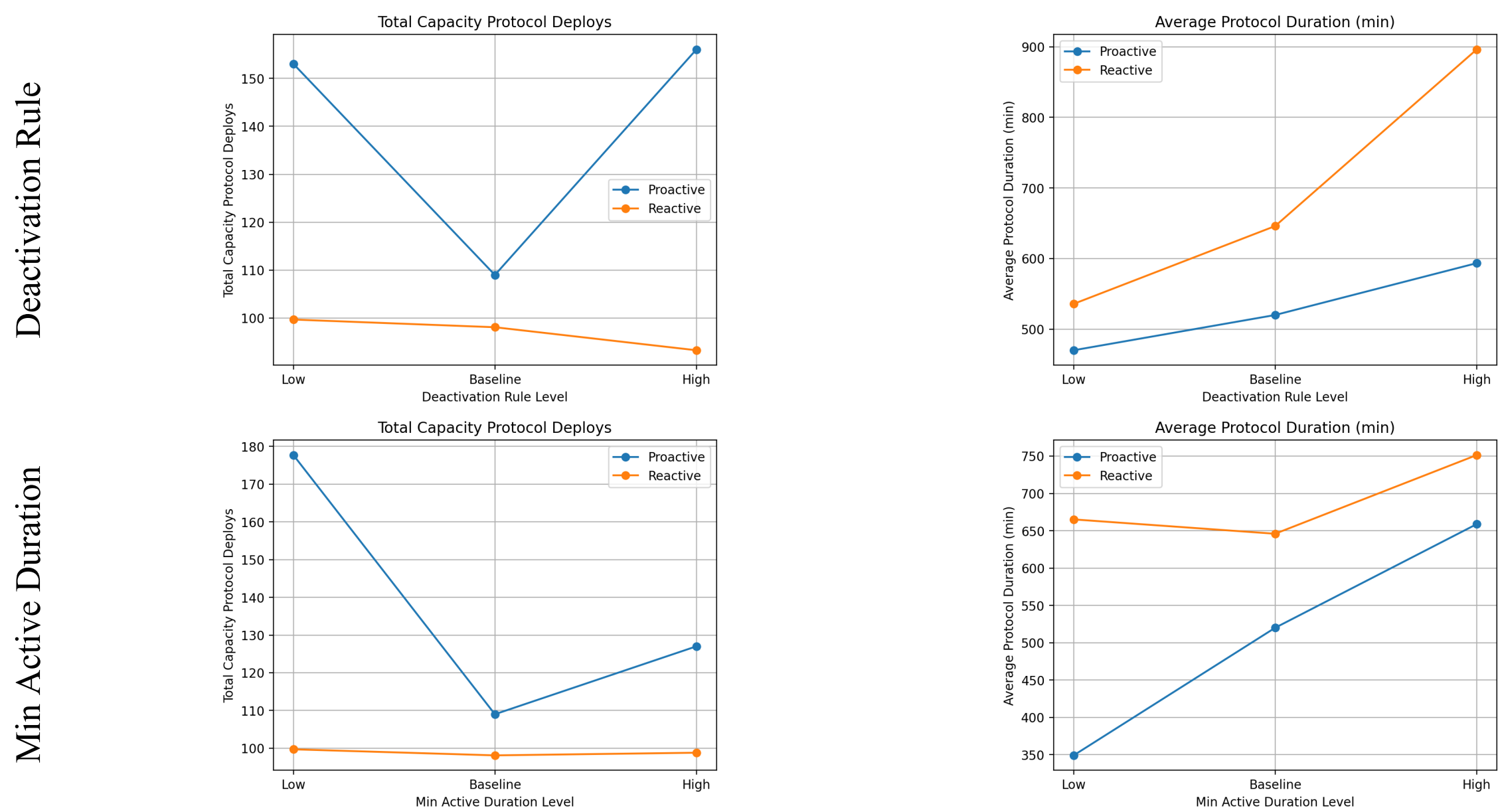

Figure 8: FCP deployment performance results from the sensitivity analysis of FCP parameters.

Figure 9 shows the ED performance results under simultaneous variation of FCP parameters, where two parameters are varied together while all other protocol parameters remain at their baseline values. This analysis shows that the performance of proactive FCP deployment depends on the joint configuration of protocol parameters. In particular, the interaction between activation delay and resource levels shows that delayed responses can reduce the benefit of proactive deployment, while higher resource levels may partially compensate for slower activation. Under high-delay conditions, proactive deployment exhibits reduced effectiveness unless additional capacity is provided. This suggests that activation timing and resource intensity jointly influence proactive FCP performance, although the magnitude and direction of this effect vary by outcome measure. Similarly, the interaction between cooldown duration and minimum active duration affects the stability of ED operations, where shorter cooldowns combined with shorter active durations enable more responsive adjustments but may increase variability if not balanced appropriately. Overall, these interaction results suggest that proactive FCP is most effective when activation timing, resource intensity, cooldown duration, and minimum active duration are calibrated together rather than optimized separately.

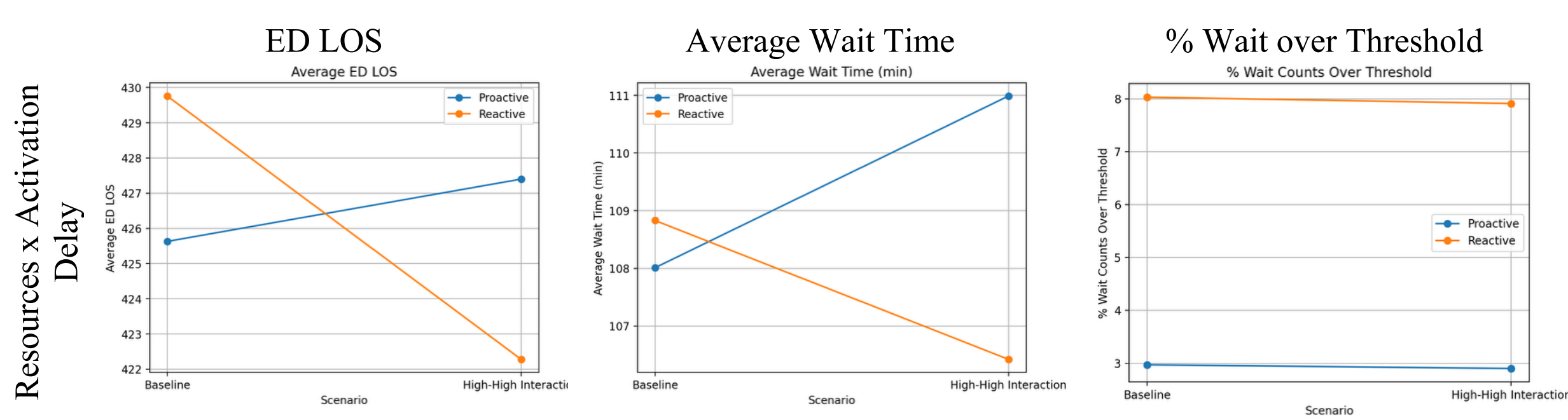

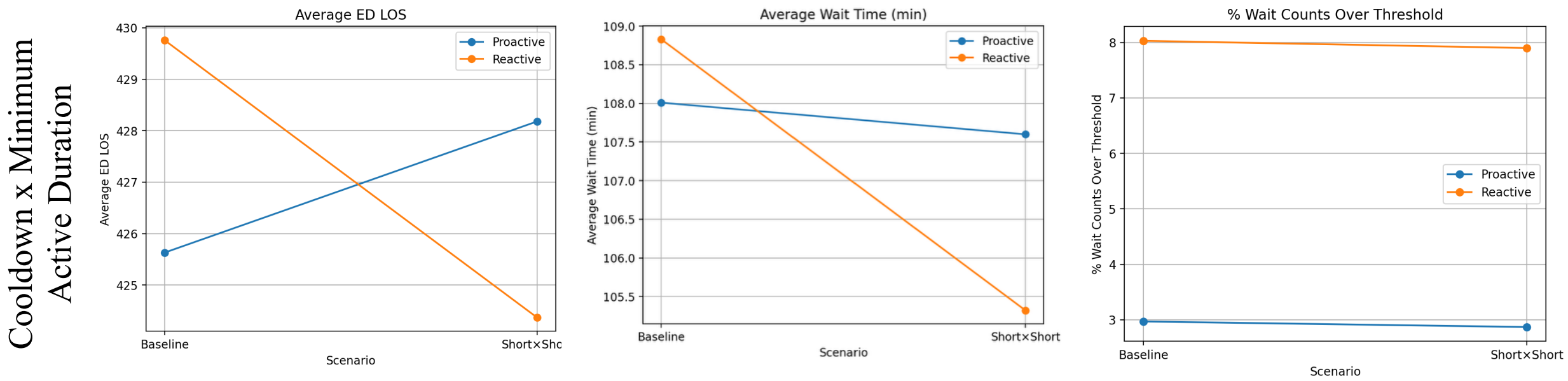

Figure 9: ED performance results for the analysis of interaction between FCP parameters.

The interaction effects on FCP deployment behavior further emphasize the importance of coordinated parameter design as shown in Figure 10. The combination of high activation delay and high resource levels results in a deployment pattern characterized by more frequent proactive activations and moderate protocol durations. This implies that timing and resource intensity can jointly influence how the protocol is used. In contrast, the combination of short cooldown and short minimum duration leads to a high-frequency, short-duration deployment pattern. These interaction-driven changes in deployment dynamics directly influence how capacity is utilized over time and provide insight into the operational trade-offs between responsiveness, resource intensity, and intervention stability.

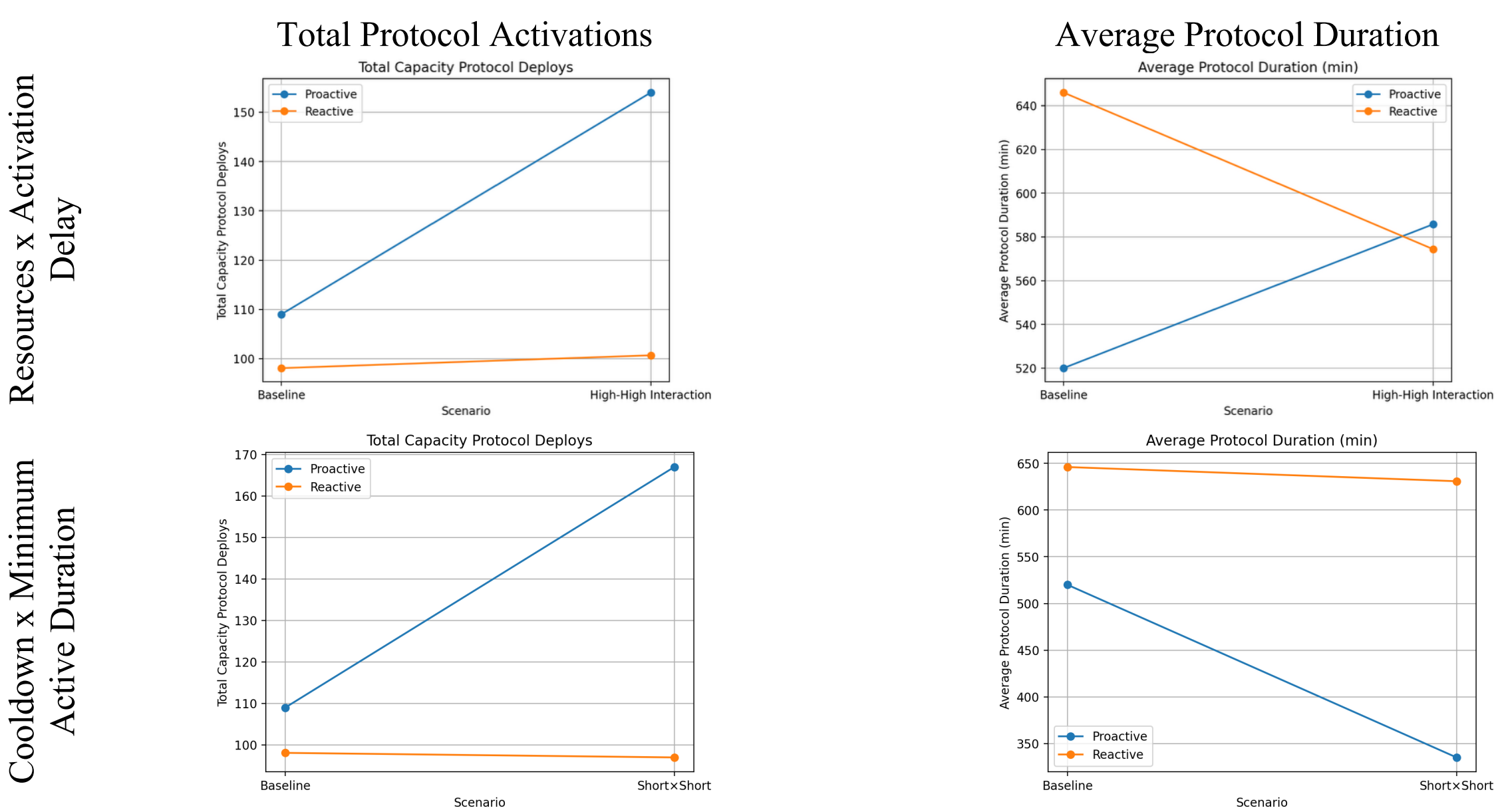

Figure 10: FCP deployment performance results for the analysis of interaction between FCP parameters.

## Discussion

In this study, we designed and evaluated a machine learning–guided decision logic for proactive ED FCP deployment using a validated simulation environment. This study is a preliminary step before the actual implementation of the proposed ML models within a decision support system to be used in the ED of a large U.S. hospital. The primary goal of this work was not only to assess the performance of proactive deployment relative to the conventional reactive approach, but also to use simulation as a safe and controlled environment to pilot-test the decision logic, operational impact, and potential risks of the system prior to real-world implementation. Because proactive deployment directly influences staffing, room availability, transfer prioritization, and inpatient coordination, any operational changes must be tested carefully. Running such experiments in a live ED is impractical and potentially disruptive. Our simulation framework therefore provides a vital intermediary step, allowing us to embed real-time machine learning predictions, model decision-making processes, and quantify system-wide effects across thousands of virtual patient encounters before piloting the tool in practice.

The simulation results demonstrate that proactive, prediction-guided FCP deployment leads to smoother and more stable ED operations across a range of crowding threshold settings. Proactive deployment consistently resulted in lower percentages of time above the crowding threshold and shorter protocol durations compared with reactive deployment. Rather than allowing crowding to escalate until a threshold is exceeded, proactive deployment intervenes earlier, reducing the likelihood that the ED remains in a prolonged overcrowded state. This is a key advantage, as past research shows that once ED crowding reaches certain thresholds, additional resources yield diminishing benefits, and crowding propagates throughout the hospital.[37]

The sensitivity analysis further refines these findings by demonstrating that the benefits of proactive FCP deployment are robust but not invariant to protocol design choices. Across a wide range of parameter configurations, proactive deployment consistently maintained better control of severe crowding exposure, as reflected by lower percentages of time above threshold. However, the magnitude of these improvements depended on parameters that control responsiveness and intervention strength. Differences in average ED LOS and average wait time were smaller and less consistent across sensitivity settings, suggesting that proactive FCP primarily improves control of extreme crowding rather than uniformly reducing all throughput measures. Parameters that influence responsiveness and persistence, including activation delay, cooldown duration, minimum active duration, resource levels, and deactivation rules, affected how frequently and how long the protocol was deployed. From an implementation perspective, successful deployment of proactive FCP requires joint tuning of decision thresholds, response delays, and resource configurations rather than isolated parameter adjustments.

By contrast, the reactive approach functions as a late-stage intervention: it responds only after crowding has already manifested. As our simulation shows, once the ED crosses these thresholds, patient queues can persist for longer periods, requiring longer protocol durations and slower recovery. In multiple scenarios, reactive deployment experienced higher percentages of time spent above crowding thresholds compared with proactive deployment. Even when FCP was deployed an equivalent number of times, the reactive strategy's interventions occurred too late to prevent the development of significant bottlenecks, resulting in more prolonged activation periods and slower system stabilization.

The operational significance of smoother ED dynamics is important for patient flow and hospital operations. Reduced variation in crowding is strongly associated with better patient experience, reduced staff burnout, improved throughput predictability, and fewer downstream delays for inpatient units.[38] Our results show that proactive deployment achieves these benefits without prolonged intervention durations. Although proactive deployment triggers FCP slightly more often under sensitive thresholds, its activations are shorter and more efficient. This suggests that earlier response may help avoid prolonged crowding cycles, although prospective clinical evaluation is needed to confirm this effect in practice.

Our findings support the feasibility and potential advantages of implementing the predictive decision support system in the partner hospital's ED. The simulation study helps identify expected performance patterns, compare deployment strategies, and evaluate whether proactive activation creates obvious operational trade-offs before live implementation. The next phase of our work will involve a controlled pilot in the live clinical environment, where the predictive system will run in real time and the decision logic will be made available to operational leaders. The simulation infrastructure developed here will serve as an ongoing testbed to evaluate alternative threshold rules, staffing responses, or prediction models prior to hospital-wide adoption.

In summary, our results indicate that proactive, prediction-guided FCP deployment offers consistent improvements in ED stability and crowding control, particularly by reducing the percentage of time spent above crowding thresholds and shortening protocol durations. Simulation served as a necessary intermediary validation step and suggested that early, data-driven interventions may reduce the escalation of crowding dynamics that often challenge EDs. These insights form the foundation for real-world implementation of this decision-support system and highlight a broader methodological pathway for integrating predictive analytics into hospital operations in a safe, evidence-based manner.

## Limitations and Future Work

This study has limitations that should be considered when interpreting the results. First, the simulation model was calibrated using operational parameters, flow patterns, and bed configurations from a single institution. Although the model reproduces that ED's dynamics in detail, the absolute performance levels and the specific FCP thresholds may not generalize to hospitals with different case mixes, staffing models, or inpatient bed configurations. Second, we evaluated proactive deployment using a single short-horizon crowding prediction model. Other prediction approaches (e.g., alternative feature sets, architectures, or retraining strategies) might yield different trigger patterns and could further enhance the benefits we observed. In addition, the proactive strategy was evaluated within a simulation environment rather than through prospective clinical deployment. Therefore, the results should be interpreted as evidence of operational feasibility and potential impact, not as definitive evidence of real-world effectiveness.

Third, the decision-support logic assumes instantaneous and full operational compliance once the FCP is activated. In practice, additional hallway beds, nursing staff, and transfer restrictions are implemented gradually, subject to staffing availability, institutional policies, and clinical judgment. Similarly, we do not explicitly model communication, coordination, or behavioral delays that might be present. These factors may reduce the effective speed and magnitude of the protocol's impact in a live environment. Fourth, we modeled FCP as a fixed intervention bundle. This allowed us to isolate the effect of activation timing, but it does not capture the full range of site-specific FCP actions, partial implementation, or clinician discretion that may occur in practice. Finally, as with all simulation studies, the model necessarily abstracts many aspects of ED operations and hospital-wide dynamics. While the model structure and parameters were chosen to reflect real processes at the study site, future work will be needed to validate the decision-support system prospectively and to test its performance in other institutional settings.

Future research may extend this work in several directions. First, rather than treating the FCP as a fixed bundle of hallway beds, additional nurses, and transfer restrictions, we plan to use the simulation model as a testbed for optimizing intervention packages. Second, the current threshold-based activation logic can be generalized using reinforcement learning, allowing an agent to learn activation and deactivation policies directly from simulated experience and adapt them to changing arrival patterns, seasonality, or hospital-wide demand. Finally, the ultimate test of the proposed decision-support system will be prospective real-world pilot studies at the partner hospital, in which the predictive model and the proactive FCP logic are integrated with existing operational dashboards and evaluated using pragmatic implementation and effectiveness endpoints, such as crowding metrics, staff acceptance, and workflow impact. Such studies will also be needed to assess alert burden, user trust, workflow fit, implementation fidelity, and unintended consequences before broader deployment.

## Acknowledgements

This study was funded by the Agency for Healthcare Research and Quality (AHRQ) under grant 5R21HS029410-02**.** The funder played no role in study design, data collection, analysis and interpretation of data, or the writing of this manuscript.

## Author Contributions

K.A. and A.A. conceptualized the study. K.A. developed the simulation model, implemented the FCP logic, conducted the experiments, generated the results, and drafted the manuscript. O.V. developed and evaluated the prediction model. J.B. provided clinical and operational expertise on ED workflow, crowding dynamics, and FCP implementation. B.O. contributed to funding acquisition as MPI and manuscript review. A.A. served as contact PI for the funded project, supervised the study, guided the methodology, interpreted the findings, and revised the manuscript. All authors reviewed and approved the final manuscript.

## Code Availability

The project code associated with this research is available on GitHub. It can be accessed through the following link: https://github.com/aramky/ED_AI_Simulation/tree/main

## Funding

**Funding: This work was supported by the Agency for Healthcare Research and Quality (AHRQ) under grant R21HS029410.**

## Conflicts of Interests

The authors declare no financial or non-financial competing interests.

## Data Availability

Data for simulation input parameters are available in the appendix.

## Supplementary Material

This appendix includes the ED data used in the simulation model and a summary of simulation parameters.

Table A1: Core simulation settings and ED operating assumptions.

| Component | Value / Rule |
|---|---|
| Simulation approach | Discrete-event simulation with event-scheduling logic |
| Simulation run length | 365 days |
| Warm-up period | 72 hours |
| Number of simulation replicates | 10 |
| ED treatment room capacity | 75 rooms |
| Triage nurse capacity | 2 nurses |
| Treatment nurse staffing schedule by 4-hour period | 7 AM–11 AM: 24<br>11 AM–3 PM: 29<br>3 PM–7 PM: 29<br>7 PM–11 PM: 30<br>11 PM–3 AM: 29<br>3 AM–7 AM: 27 |
| Nurse-to-patient ratio by ESI level | ESI 1: 1, ESI 2: 1:2, ESI 3: 1:3, ESI 4: 1:4, ESI 5: 1:4 |
| Number of inpatient units | 28 |
| FCP hallway beds added | 10 |
| FCP temporary nurses added | 5 |
| Transfer cancellation under FCP | 3% of patient arrivals |
| Proactive prediction frequency | Every 60 minutes |
| Prediction horizon | 6 hours ahead |

Table A2: Time-varying Poisson ED arrival rates.

| Clock time | Arrival rate (patients/minute) |
|---|---|
| 7 AM–11 AM | 0.1335 |
| 11 AM–3 PM | 0.1825 |
| 3 PM–7 PM | 0.1572 |
| 7 PM–11 PM | 0.1209 |
| 11 PM–3 AM | 0.1058 |
| 3 AM–7 AM | 0.0673 |

Table A3: ED front-end process times and patient acuity distributions.

| Input element | Distribution / Rule | Parameters / Probabilities |
|---|---|---|
| Registration time | Uniform | [5, 7] min |
| Triage time | Uniform | [10, 15] min |
| ESI distribution | Discrete Empirical | $P(ESI) = \begin{cases} 0.0226, & ESI\ 1 \\ 0.3214, & ESI\ 2 \\ 0.4870, & ESI\ 3 \\ 0.1458, & ESI\ 4 \\ 0.0232, & ESI\ 5 \end{cases}$ |
| Admission probability by ESI | Bernoulli by ESI | ESI 1: 0.848<br>ESI 2: 0.505<br>ESI 3: 0.195<br>ESI 4–5: 0 |

Table A4: ED treatment-time distributions by ESI level (in minutes).

| ESI level | Distribution | Parameters |
|---|---|---|
| ESI 1 | Beta | $a = 1.8180$, $b = 5.5030$, $loc = 14.3658$, $scale = 556.3718$ |
| ESI 2 | Beta | $a = 1.3392$, $b = 1.3782$, $loc = 14.8682$, $scale = 391.6886$ |
| ESI 3 | Beta | $a = 1.4159$, $b = 2.9178$, $loc = 14.2525$, $scale = 687.7421$ |
| ESI 4 | Gamma | $a = 1.4501$, $loc = 14.1376$, $scale = 73.8015$ |
| ESI 5 | Exponential | $loc = 15.0$, $scale = 67.4706$ |

Table A5: Inpatient units and bed capacities.

| Unit | Beds | Unit | Beds |
|---|---|---|---|
| Unit 1 | 29 | Unit 15 | 24 |
| Unit 2 | 21 | Unit 16 | 26 |
| Unit 3 | 21 | Unit 17 | 18 |
| Unit 4 | 32 | Unit 18 | 25 |
| Unit 5 | 29 | Unit 19 | 36 |
| Unit 6 | 31 | Unit 20 | 20 |
| Unit 7 | 31 | Unit 21 | 21 |
| Unit 8 | 28 | Unit 22 | 12 |
| Unit 9 | 30 | Unit 23 | 19 |
| Unit 10 | 30 | Unit 24 | 28 |
| Unit 11 | 31 | Unit 25 | 28 |
| Unit 12 | 31 | Unit 26 | 21 |
| Unit 13 | 20 | Unit 27 | 21 |
| Unit 14 | 32 | Unit 28 | 21 |

Table A6: Inpatient unit assignment probabilities by ESI level

| Unit | ESI 1 | ESI 2 | ESI 3 |
|---|---|---|---|
| Unit 1 | 0.0067 | 0.0153 | 0.0042 |
| Unit 2 | 0.0450 | 0.0124 | 0.0055 |
| Unit 3 | 0.0492 | 0.0128 | 0.0042 |
| Unit 4 | 0.0308 | 0.0153 | 0.0133 |
| Unit 5 | 0.0125 | 0.0185 | 0.0212 |
| Unit 6 | 0.0034 | 0.0169 | 0.0274 |
| Unit 7 | 0.0850 | 0.0209 | 0.0027 |
| Unit 8 | 0.0858 | 0.0190 | 0.0084 |
| Unit 9 | 0.0375 | 0.0292 | 0.0086 |
| Unit 10 | 0.0050 | 0.0243 | 0.0351 |
| Unit 11 | 0.0092 | 0.0214 | 0.0411 |
| Unit 12 | 0.0050 | 0.0199 | 0.0462 |
| Unit 13 | 0.0725 | 0.0272 | 0.0159 |
| Unit 14 | 0.0158 | 0.0224 | 0.0446 |
| Unit 15 | 0.0200 | 0.0279 | 0.0345 |
| Unit 16 | 0.0092 | 0.0280 | 0.0623 |
| Unit 17 | 0.0100 | 0.0299 | 0.0583 |
| Unit 18 | 0.0242 | 0.0308 | 0.0556 |
| Unit 19 | 0.0317 | 0.0497 | 0.0205 |
| Unit 20 | 0.0167 | 0.0367 | 0.0590 |
| Unit 21 | 0.0134 | 0.0494 | 0.0369 |
| Unit 22 | 0.0267 | 0.0430 | 0.0738 |
| Unit 23 | 0.0217 | 0.0533 | 0.0552 |
| Unit 24 | 0.0433 | 0.0471 | 0.0698 |
| Unit 25 | 0.0367 | 0.0610 | 0.0548 |
| Unit 26 | 0.0808 | 0.0940 | 0.0243 |
| Unit 27 | 0.0700 | 0.0632 | 0.0961 |
| Unit 28 | 0.1324 | 0.1102 | 0.0208 |

**Note:** ESI 4–5 patients are assumed to have zero admission probabilities.

Table A7: Mean of exponentially distributed inpatient length of stay (minutes) by unit and ESI level.

| Unit | ESI 1 | ESI 2 | ESI 3 |
|---|---|---|---|
| Unit 1 | 103356.00 | 8472.14 | 8244.57 |
| Unit 2 | 26362.68 | 9126.09 | 6624.65 |
| Unit 3 | 5789.26 | 7790.09 | 12701.52 |
| Unit 4 | 21751.38 | 10359.51 | 8803.42 |
| Unit 5 | 14394.00 | 13722.38 | 9109.80 |
| Unit 6 | 18402.25 | 10469.60 | 7533.92 |
| Unit 7 | 38037.67 | 8685.92 | 7221.49 |
| Unit 8 | 9538.12 | 9569.57 | 6214.50 |
| Unit 9 | 31291.00 | 17381.36 | 11159.79 |
| Unit 10 | 30434.83 | 12429.59 | 10902.08 |
| Unit 11 | 21566.93 | 11068.85 | 9457.68 |
| Unit 12 | 17450.94 | 12625.39 | 11265.07 |
| Unit 13 | 7806.40 | 18596.85 | 7285.09 |
| Unit 14 | 14859.15 | 8100.26 | 11730.62 |
| Unit 15 | 11981.38 | 12962.04 | 12837.20 |
| Unit 16 | 18200.48 | 12258.58 | 9424.27 |
| Unit 17 | 5912.53 | 8844.78 | 18251.18 |
| Unit 18 | 7081.51 | 11293.35 | 17688.47 |
| Unit 19 | 9237.84 | 8437.51 | 4587.91 |
| Unit 20 | 15764.74 | 7474.61 | 8560.51 |
| Unit 21 | 20516.35 | 12390.81 | 9902.58 |
| Unit 22 | 20681.50 | 20179.12 | 16122.08 |
| Unit 23 | 53776.76 | 20091.99 | 12284.73 |
| Unit 24 | 11877.44 | 10347.87 | 5453.25 |
| Unit 25 | 31018.28 | 10060.82 | 6792.97 |
| Unit 26 | 41656.54 | 13682.45 | 8048.06 |
| Unit 27 | 37480.58 | 13674.18 | 12587.05 |
| Unit 28 | 15889.82 | 9613.30 | 8713.34 |

Table A8: Prediction model inputs and FCP control rules.

| Item | Details |
|---|---|
| Prediction model | SVM model loaded from pickle file: (`wait_predictor_model_SVM.pkl`). |
| Prediction inputs | Year, month, day of month, hour, day of week, extreme-time flag, and 24 lagged hourly wait counts. |
| Extreme-time flag | 1 if most recent hourly wait count > $\mu + 2\ \sigma$ of logged wait counts; 0 otherwise. |
| Proactive trigger rule 1 | Predicted crowding exceeds threshold for 2 consecutive hours |
| Proactive trigger rule 2 | Increasing wait-count trend over 4 hours and treatment room utilization > 80%. |
| Reactive trigger rule | Current wait count ≥ threshold. |
| FCP cooldown | 6 hours. |
| Minimum active duration | 6 hours. |
| Deactivation rule | Wait count below 60% of threshold for 4 consecutive hours |
| Activation delay | 3 hours. |